\documentstyle[12pt,amssymb,amscd]{amsart} 
\author{Dmitry E. Tamarkin}
\title{Formality of Chain Operad of Small Squares}
\newtheorem{Theorem}{THEOREM}[section]

 \newtheorem{Proposition}[Theorem]{PROPOSITION}

\def \g {{\frak g}}
\begin{document}
\hsize 15cm
\vsize 20cm
\leftmargin -2cm
\topmargin=1cm
\maketitle
\def \cs {C^{\rm sing}_\bullet}
\def \Q{{\Bbb Q}}
\def \dash{\mbox{-}}
\input amssym.def
\section{Introduction} We prove that the chain  operad
of small squares is formal. Together with the Deligne
conjecture about the action 
of this operad on Hochschild cochains of an
 associative algebra, this implies existence of a structure of a 
homotopy Gerstenhaber algebra on Hochschild cochains. 
This fact clarifies  situation
with the proof of M. Kontsevich formality theorem in the 
paper of the author \ref{Tam}. 
The formality
of the operad  follows quite easily from the existence of an associator.
The author would like to thank 
Boris Tsygan, Paul Bressler, and Maxim Kontsevich fior their help.

{\bf Remark}  Maxim Kontsevich
 found a proof of a more general result:
chain operad of small balls is forlmal in all dimensions as an operad
of coalgebras. Also, he has pointed out that the construction presented here
allowes one to show the formality on the level of coalgebras as well,
since all the maps of operads involved are the maps of operads of coalgebras.

\section{Small Square Operad (after \cite{F})} We reproduce the construction
of the small square operad from \cite{F}. First, the symmetric groups in the 
definition of operad are replaced with the braid groups and we obtain the
 notion  of braided operad. A {\em topological $B_\infty$-operad} 
$X$ is defined as a braided operad such that all its spaces $X(n)$
are contractible and the braid group $B_n$ acts freely on $X(n)$.
If $X$ and $Y$ are topological $B_\infty$-operads, then so is
$X\times Y$ and we have homotopy equivalences 
\begin{equation}\label{equiv}
p_1:X\times Y\to X;\quad p_2:X\times Y\to Y,
\end{equation}
where $p_1,p_2$ are the projections.

Let $PB_n$ be the group of pure braids with $n$ strands.
Given a  topological $B_\infty$-operad $X$,
the corresponding {\em operad of small squares} 
is a symmetric operad $X'$ such that $X'(n)=X(n)/PB_n$ with the induced structure maps.  The maps (\ref{equiv}) guarantee that any two operads of small
squares are connected by a chain of homotopy equivalences. It is proven in \cite{F} that the classical operad of May (whose $n-th$ space is the configuration space of $n$ disjoint numbered squares inside the unit square such that the corresponding sides are parallel) is an operad of small squares in our sense. 

The functor of singular chains $C^{\rm sing}_\bullet:Top\to Complexes$ has a natural
tensor structure  given by the Eilenberg-Zilber map $EZ:\cs(X)\otimes \cs(Y)\to 
\cs(X\times Y)$. Therefore, for a topological operad $O$,
 the collection $\cs(O(\bullet))$ has a structure of a $dg$-operad.
The structure map of the $i$-th
insertion is 
$$
\cs(O(n))\otimes \cs(O(m))\stackrel{EZ}{\to} \cs(O(n)\times O(m))\stackrel{o_{i*}^O}{\to}\cs(O(n+m-1)),
$$
where $o_i^O$ is the structure map of the $i$-th insertion in $O$.

 For a small square operad $X$ consider the operad $E_2(X)=\cs(X)$. Any two such operads are quasi-isomorphic, where quasi-isomorphic
means connected by a chain of quasi-isomorphisms.
In particular, the homology operad of any of $E_2(X)$
is the operad $e_2$ controlling Gerstenhaber algebras
(see section \ref{apbn} for the definition of $e_2$). Our goal is to show that 
\begin{Theorem}
Any operad $E_2(X)$ is quasi-isomorphic to its homology  operad $e_2$.
\end{Theorem}  
For this it suffices to pick a particular 
small square operad $X$.
\section{Realization of $E_2$}
\subsection{Operad of categories $PaB_n$} (after \cite{BN})
First, let us  reproduce from \cite{BN} the definition of the category $PaB_n$.
Let $B_n$ $(PB_n)$ be the group
of braids (pure braids) with $n$ strands, let $S_n$ be the symmetric group.
 Let $p:B_n\to S_n$
be the canonical projection with the kernel $PB_n$. We assume that the strands
of any braid are numbered in the order determined by their origins.

The objects of the category $PB_n$ are 
 parenthesized permutations of elements $1,2,\ldots,n$ (that is pairs $(\sigma, p)$,
where $\sigma\in S_n$ and $\pi$ is a parenthesation of the non-associative product of $n$ elements). The morphisms between $(\sigma_1,\pi_1)$ and
$(\sigma_2,\pi_2)$ are such braids from $B_n$ that any strand joints an element
of $\sigma_1$ with the same element of $\sigma_2$, in other words,
${\rm Mor}((\sigma_1,\pi_1),(\sigma_2,\pi_2))=p^{-1}(\sigma_2^{-1}\sigma_1)$. 
The composition law is induced from the one on $B_n$.
 The symmetric group
$S_n$ acts on $PaB_n$ via renumbering the objects $T_\sigma(\sigma_1,\pi_1)=
(\sigma\sigma_1,\pi_1)$ and acts identically on morphisms.

The collection of categories $PaB_n$ form an operad.
Indeed, the collection ${\rm Ob} PaB_\bullet$ forms a free operad
in the category of sets generated by one binary noncommutative operation. Let us describe the structure map
$o_k$ of the insertion into the $k$-th position on the level of morphisms.
Suppose we insert $y:(\sigma_1,\pi_1)\to (\sigma_2,\pi_2)$ into $x:(\sigma_3,\pi_3)\to (\sigma_4,\pi_4)$. We replace the strand number $\sigma^{-1}_3(k)$ of the braid $x$ by the brade
$y$ made very narrow.
\subsection{Operad of classifying spaces} We have the functor of taking the nerve
$N:{\em Cat}\to \Delta^oEns$ and the functor of topological realization
$|\ |:\Delta^oEns\to Cellular Spaces $. These functors  behave well with respect to the 
symmetric monoidal structures, therefore   the collection of cellular 
complexes $X_n=|NPaB_n|$ forms a cellular operad. 
One checks that this operad is a small square operad. Indeed, let $PaB_n'$
be the category whose objects are pairs $(x,y)$, where $x$ belongs to the braid group $B_n$ and
$y$ is a parenthesation of the non-associative product of $n$ elements, and 
there is a unique morphism between any two objects. We have  a free left action
of $B_n$ on $PaB_n'$: $(x,y)\to (gx,y))$ and a 
braided operad structure on $PaB_\bullet'$ 
(the structure maps are defined similarly to $PaB_n$).  One checks that the corresponding operad of classifying spaces is a 
 topological $B_\infty$-operad and that the corresponding small square operad
is isomorphic to $X_\bullet$. 
 
Consider the corresponding 
chain operad. Let $C_\bullet(NPaB_n)$ 
be the chain complex over $\Bbb Q$ of $NPaB_n$
as a simplicial set. The collection $C_\bullet(NPaB_n)$ forms a $dg$-operad 
(via the Eilenberg-Zilber 
map). Since $C_\bullet(NPaB_n)$ is just a bar complex
of the category $PaB_n$,
this operad will be denoted by $C_\bullet(PaB_\bullet)$. We have a canonical quasi-isomorphism of operads 
$C_\bullet(PaB_\bullet)\to \cs|NPaB_\bullet|$. Therefore, it suffices to construct a 
quasi-isomorphism of $C_\bullet(PaB_\bullet)$ and $e_2$.

\section{Operad of algebras $A^{pb}_n$ and construction of quasi-isomorphism}
\label{apbn}
 By definition \cite{Dr} $A^{pb}_n$ is the algebra over $\Bbb Q$ of
power series in the noncommutative  variables 
\begin{equation}\label{tij} 
t_{ij},1\leq i,j\leq n;\quad i\neq j; \quad t_{ij}=t_{ji}
\end{equation}
with relations 
\begin{equation}\label{rel}
[t_{ij}+t_{ik},t_{jk}]=0.
\end{equation} 
Let $I_n$ be the double-sided ideal generated by all $t_{ij}$.
We have a canonical projection 
\begin{equation}\label{chi}
\chi: A^{pb}_n\to A^{pb}_n/I_n\cong {\Bbb Q}.
\end{equation}
The symmetric group $S_n$
acts naturally on $A^{pb}_n$ so that $T_\sigma t_{ij}=t_{\sigma(i)\sigma(j)}$.
The collection $A^{pb}_n$ forms an operad in the category of algebras
in a well-known way. The map of the insertion into the $i$-th
position $o_i:A^{pb}_n\otimes A^{pb}_m\to A^{pb}_{n+m-1}$ looks as follows.

Let $$
\phi(k)=\left\{\begin{matrix} k, && k\leq i;\\
                              k+m-1,&& k>i.
               \end{matrix}\right.
    $$
Then                        
$$
o_i(t_{pq}\otimes 1)=\left\{\begin{matrix} t_{\phi(p)\phi(q)}, &&  p,q\neq i;\\
                                           \sum_{r=i}^{i+m-1} t_{r\phi(q)}, && p=i;
                            \end{matrix}\right. 
$$   
$$
o_i(1\otimes t_{pq})=t_{i+p-1,i+q-1}.
$$

Any algebra with unit over $\Bbb Q$ gives rise to a $\Bbb Q$-additive category
$C_A$ with one object. Denote by $\Q Cat$ the category of small $\Q$-additive categories, and by
 $\Q Cat'=\Q Cat/C_\Q$ the over-category of $\Q Cat$  over $C_\Q$. Its objects are the elements of ${\rm Mor}_{\Q Cat}(x,C_\Q)$, where $x\in \Q Cat$.
 A morphism between $\phi$ and $\psi$, where
$\phi:x\to C_\Q$; $\psi:y\to C_\Q$, is a morphism $\sigma:x\to y$ in 
$\Q Cat$ such that
$\sigma\psi=\phi$. This category has a clear symmetric monoidal structure.
 We have the functor of nerve $N^{\Q}:\Q Cat'\to \Delta^oVect$, which is
the straight analogue of the nerve of an arbitrary category,
and the functor $C_\bullet:\Delta^oVect\to Complexes$. Both of these functors have tensor structure (on the latter functor it is defined
via the Eilenberg-Zilber map), therefore we have a through functor
$\Q Cat'\to Compexes$ and the induced functor 
$$\Q Cat'\dash Operads\to dg\dash Operads,
$$
which will be denoted by $C_\bullet^{\Q}$.

The map (\ref{chi}) produces a morphism $\chi_*:C_{A^{pb}_n}\to C_\Q$ and defines an object $O_A(n)\in \Q Cat'$. The operad structure on
$A^{pb}_n$ defines an operad structure on the collection $O_A(n)$.
The complex $C_\bullet^{\Q}(O_A(n)) $ looks as follows:
 $C_n^{\Q}(O_A(k))\cong A_k^{pb \otimes n}$;  
$$
d a_1\otimes\cdots\otimes a_n=\chi(a_1)a_2\otimes\ldots\otimes a_n-
a_1a_2\otimes\cdots\otimes a_n+\ldots+(-1)^{n-1}a_1\otimes\ldots\otimes a_{n-1}
\chi(a_n).
$$ 
This is the bar complex for
${\rm Tor}^{A^{pb}_n}(A^{pb}_n/I_n,A^{pb}_n/I_n) $.

Let $e_2$ be the operad of graded vector spaces governing the Gerstenhaber
algebras. It is generated by two binary operations:
the commutative associative multiplication of degree zero,
 which is denoted by $\cdot$, and the commutative bracket of degree -1
denoted by $\{,\}$.

These operations satisfy the Leibnitz identity $\{ab,c\}=a\{b,c\}+(-1)^{b(c+1)}
\{a,c\}b$ and the Jacoby identity 
$$
(-1)^{|a|}\{a,\{b,c\}\}+(-1)^{|a||b|+|b|}\{b,\{a,c\}\}+
(-1)^{|a||c|+|b||c|+|c|}\{c,\{a,b\}\}=0.
$$
  We have a morphism of operads
\begin{equation}\label{k}
k:e_2\to C_\bullet^{\Q}(O_A),
\end{equation}
 which is defined on $e_2(2)$
as follows: $k(\cdot)=1\in C_0^{\Q}(O_A)$; $k(\{,\})=t_{12}\in C_1^{\Q}(O_A)$. 
Direct check shows that this map respects the relations in $e_2$.
\begin{Proposition}\label{qis}
The map $k$ is a quasi-isomorphism of operads
\end{Proposition}

\noindent{\em Proof}. Let ${\g}_n$ be the graded Lie algebra
generated by  the elements (\ref{tij}) and relations (\ref{rel}), and
the grading is defined by setting $|t_{ij}|=1$. Then the
universal enveloping algebra $U\g_n$ is a graded associative algebra, and
$A^{pb}_n$ is the completion of $U\g_n$ with respect to the grading.
The algebras $U\g_\bullet$ form an operad with the same structure maps as in
$A^{pb}_\bullet$. The inclusion 
\begin{equation}\label{incl}
U\g_n\to A^{pb}_n
\end{equation} is a morphism
of operads. We have a canonical projection $\chi:U\g_n\to \Bbb Q$, therefore
the copllection $C_{U\g_\bullet}$
forms an operad in $\Q Cat'$ and we have a dg-operad  $C_\bullet^{\Q} C_{U\g_\bullet}$ which will be denoted by $C_\bullet^{\Q} U\g_\bullet$. The injection (\ref{incl})
induces a morphism of operads 
\begin{equation}\label{incl*}
i_*:C_\bullet^{\Q}(U\g_\bullet)\to C_\bullet^{\Q}(O_A).
\end{equation}

It is clear that ${\rm Tor}^{A^{pb}_n}_\bullet(A^{pb}_n/I_n,A^{pb}_n/I_n)$
is the same as the completion of $H_\bullet(\g_n)\cong {\rm Tor}^{U\g_n}_\bullet({\Bbb Q},{\Bbb Q})\cong H_\bullet(C_\bullet^{\Q}(U\g_n))$ with respect to the grading induced from $\g_n$. 

We have a natural injection $\g_{n-1}\to \g_n$. One sees that
the Lie subalgebra $\i_n\subset \g_n$ generated by
$t_{nk}$, $k=1,\ldots n-1$ is free and is an ideal in $\g_n$.
Also, we have $\g_n=\i_n\oplus \g_{n-1}$
in the category of vector spaces. The Serre-Hochschild spectral sequence
$E^2_{\bullet,\bullet}=H_\bullet(\g_{n-1},H_\bullet \i_n)\Rightarrow H_\bullet(\g_n)$ collapses
at $E^2$ and shows that 
\begin{equation}\label{decompose}
H_\bullet(\g_n)\cong H_\bullet(\g_{n-1})\oplus
(\oplus_{k=1}^{n-1}H_\bullet(\g_{n-1}))[-1],
\end{equation} where the first summand is the image of
$H_\bullet(\g_{n-1})$
 under the injection $\g_{n-1}\to \g_n$. Let us describe
the remaining $n-1$ summands.

 Note that $\g_2$
is one-dimensional, therefore
$H_0(\g_2)=\Bbb Q;\ H_1(\g_2)=\Bbb Q[-1]$;
$H_i(\g_2)=0$ for $i>1$.
The induction shows that the homology
of $\g_n$ is finite dimensional, therefore {\em the map (\ref{incl*})
is a quasi-isomorphism.}
The operadic maps of  insertion into
the $k$-th position $o_k:U\g_{n-1}\otimes U\g_2\to U\g_n$, where
$k=1,2,\ldots,n-1$ induce maps $o_k^*:H_\bullet(\g_{n-1})\otimes H_\bullet(\g_2)\to H_\bullet(\g_n)$, and the $(k+1)$-th summand in
(\ref{decompose})
is equal to  $o_k^*(H_\bullet(\g_{n-1})\otimes H_1(\g_2))$. The induction 
argument shows that 
\begin{enumerate}
\item[1.] The homology operad 
$n\mapsto H_\bullet(\g_n)\cong  H_\bullet(C_\bullet^{\Q}(U\g_n))$ is generated by $H_\bullet(\g_2)$,
therefore the homology operad of $C_\bullet^{\Q}(O_A(n))$
is generated by the homology of $C_\bullet^{\Q}(O_A(2))$.
\item[2.] The total dimension of $H_\bullet(\g_n)$ and of the homology
of $C_\bullet^{\Q}(O_A(n))$ is $n!$.
\end{enumerate}
The first statement implies that the map (\ref{k})
is surjective
on the homology level, and the second statement means that the map (\ref{k})
is bijective since ${\rm dim}\; e_2(n)=n!$. $\bigtriangleup$

Let ${\Bbb Q}(PaB_n)$ be the $\Q$-additive category generated by $PaB_n$. 
 We have a map ${\Bbb Q}(PaB_n)\to C_\Q$ sending all morphisms from
$PaB_n$ to ${\rm Id}$. Thus, $\Q(PaB_n)\in \Q Cat'$. 
The oreadic sructure on $PaB_\bullet$ induces the one on $\Q PaB_\bullet$.

 Any associator $\Phi\in A^{pb}_3$ over $\Bbb Q$ produces a map of operads
$\phi:{\Bbb Q}(PaB_\bullet)\to O_A(\bullet)$. 
Indeed, define $\phi$ on ${\rm Ob}\; PaB_n$ by sending any object to
the only object of $O_A(n)$.
There are only two objects in $PaB_2$, let us denote them
$x_1x_2$ and $x_2x_1$. The morphisms between these two objects correspond
to the  non-pure braids. Let $x\in B_2$ be the generator.
We define $\phi(x)=e^{t_{12}/2}$. Take the two objects $(x_1x_2)x_3$
and $x_1(x_2x_3)$ of $PaB_3$ corresponding to the identical permutation
$e\in S_3$, and the morphism $i$ between them, corresponding to the
identical braid $e_b\in B_3$. Define $\phi(i)=\Phi$. Since the operad
$PaB_\bullet$ is generated  by $x$ and $i$, these conditions define $\phi$
uniquely. The definition of the  associator is equivalent to the fact
that $\phi$ is well-defined. This construction is very similar to
the one from \cite{BN}.    

The map  $\phi$ produces a map of
operads  $C_\bullet^{\Q}{PaB_\bullet}\cong C_\bullet^{\Q}(\Q(PaB_\bullet))\to C_\bullet(O_A)$.
It is well known that the homology operad
of $C_\bullet{PaB_\bullet}$ is $e_2$.
 It is easy to check that $\phi$ is a quasi-isomorphism for $\bullet=2$ and hence
it is a quasi-isomorphism of operads
(since $e_2$ is generated by $e_2(2)$).  By Proposition 
\ref{qis}, $k$ is a quasi-isomorphism. Thus, the chain operad ${C_\bullet{PaB_\bullet}}$ is quasi-isomorphic to $e_2$.

\end{document}